\documentclass[12pt,reqno,a4paper]{amsart}
\usepackage{amsfonts,latexsym,amscd,amsxtra}
\usepackage{amsmath,amssymb}
\usepackage{ifthen}
\usepackage{graphicx,bm}
\usepackage[all]{xy}
\usepackage{color}

\nonstopmode \numberwithin{equation}{section}
\makeatletter
\newcommand*\bigcdot{\mathpalette\bigcdot@{.5}}
\newcommand*\bigcdot@[2]{\mathbin{\vcenter{\hbox{\scalebox{#2}{$\m@th#1\bullet$}}}}}
\makeatother

\newtheorem{thm}{Theorem}[section]
\newtheorem{lem}{Lemma}[section]
\newtheorem{cor}{Corollary}[section]
\newtheorem{prop}{Proposition}[section]

\newtheorem{cl}{Claim}
\newtheorem{ca}{Case}
\newtheorem{sca}{Case}
\newtheorem{scl}{Subclaim}
\newtheorem{conj}{Conjecture}

\theoremstyle{definition}
\newtheorem{defn}{Definition}

\newtheorem{op}[equation]{Open Problem}
\newtheorem{ques}[equation]{Question}
\newtheorem{rem}{Remark}[section]
\newtheorem{exam}[equation]{Example}

\newcounter {own}
\def\theown {\thesection       .\arabic{own}}

\newenvironment{pf}[1][]{%
	\vskip 3mm
	\noindent
	\ifthenelse{\equal{#1}{}}%
	{{\slshape Proof. }}%
	{{\slshape #1.} }%
}%
{\qed\bigskip}

\newcounter{alphabet}

\newenvironment{Thm}[1][]{\refstepcounter{alphabet}%
	\noindent%
	{\bf Theorem \Alph{alphabet}}%
	\ifthenelse{\equal{#1}{}}{}{ (#1)}%
	{\bf .} \itshape}{\vskip 8pt}


\newcommand{\p}{{\alpha}}
\newcommand{\q}{{\beta}}

\newcommand{\C}{{\mathbb C}}
\newcommand{\D}{{\mathbb D}}
\newcommand{\T}{{\mathbb T}}
\newcommand{\R}{{\mathbb R}}

\newcommand{\Z}{{\mathbb Z}}

\def\be{\begin{equation}}
\def\ee{\end{equation}}

\newcommand{\bee}{\begin{enumerate}}
	\newcommand{\eee}{\end{enumerate}}

\newcommand{\blem}{\begin{lem}}
	\newcommand{\elem}{\end{lem}}
\newcommand{\bthm}{\begin{thm}}
	\newcommand{\ethm}{\end{thm}}
\newcommand{\bcor}{\begin{cor}}
	\newcommand{\ecor}{\end{cor}}
\newcommand{\beg}{\begin{exam}}
	\newcommand{\eeg}{\end{exam}}
\newcommand{\begs}{\begin{examples}}
	\newcommand{\eegs}{\end{examples}}
\newcommand{\bdefe}{\begin{defn}}
	\newcommand{\edefe}{\end{defn}}
\newcommand{\bprob}{\begin{prob}}
	\newcommand{\eprob}{\end{prob}}
\newcommand{\bques}{\begin{ques}}
	\newcommand{\eques}{\end{ques}}
\newcommand{\bei}{\begin{itemize}}
	\newcommand{\eei}{\end{itemize}}
\newcommand{\bcon}{\begin{conj}}
	\newcommand{\econ}{\end{conj}}
\newcommand{\bop}{\begin{op}}
	\newcommand{\eop}{\end{op}}

\newcommand{\bca}{\begin{ca}}
	\newcommand{\eca}{\end{ca}}
\newcommand{\bsca}{\begin{sca}}
	\newcommand{\esca}{\end{sca}}

\newcommand{\bcl}{\begin{cl}}
	\newcommand{\ecl}{\end{cl}}

\newcommand{\bscl}{\begin{scl}}
	\newcommand{\escl}{\end{scl}}

\newcommand{\bcons}{\begin{conjs}}
	\newcommand{\econs}{\end{conjs}}
\newcommand{\bprop}{\begin{propo}}
	\newcommand{\eprop}{\end{propo}}
\newcommand{\br}{\begin{rem}}
	\newcommand{\er}{\end{rem}}
\newcommand{\brs}{\begin{rems}}
	\newcommand{\ers}{\end{rems}}
\newcommand{\bo}{\begin{obser}}
	\newcommand{\eo}{\end{obser}}
\newcommand{\bos}{\begin{obsers}}
	\newcommand{\eos}{\end{obsers}}
\newcommand{\bpf}{\begin{pf}}
	\newcommand{\epf}{\end{pf}}
\newcommand{\ba}{\begin{array}}
	\newcommand{\ea}{\end{array}}
\newcommand{\beq}{\begin{eqnarray}}
\newcommand{\beqq}{\begin{eqnarray*}}
\newcommand{\eeq}{\end{eqnarray}}
\newcommand{\eeqq}{\end{eqnarray*}}

\newcommand{\ds}{\displaystyle}

\begin{document}
\bibliographystyle{amsplain}
\author{Adel Khalfallah}
\address{Department of Mathematics, King Fahd University of Petroleum and
	Minerals, Dhahran 31261, Saudi Arabia}\email{khelifa@kfupm.edu.sa}
	\author{Mohamed Mhamdi}
\address{Ecole supérieure des Sciences et  de la Technologie de Hammam Sousse (ESSTHS)\\Université de Sousse}
\email{mhamdimed7@gmail.com}

\title [] {The  first  partial derivatives  of  generalized harmonic  functions}
\maketitle

\def\thefootnote{}
\footnotetext{ \texttt{\tiny
		Date: \number\day-\number\month-\number\year. }}
\makeatletter\def\thefootnote{\@arabic\c@footnote}\makeatother
\begin{abstract}
Suppose $\p,\q \in \R\backslash \Z^-$ such that  $\p+\q>-1$ and $1\leq p \leq \infty$. Let $u=P_{\p,\q}[f]$ be an $(\p,\q)$-harmonic mapping on $\D$, the unit disc of $\C$, with the boundary $f$ being absolutely continuous  and $\dot{f}\in L^p(0,2\pi)$, where $\dot{f}(e^{i\theta}):=\frac{d}{d\theta}f(e^{i\theta})$. In this paper, we investigate the membership of the partial derivatives $\partial_z u$ and $\partial_{\overline{z}}u$ in the space $H_G^{p}(\mathbb{D})$, the generalized Hardy space. We prove, if $\p+\q>0$, then both $\partial_z u$ and $\partial_{\overline{z}}u$  are in $H_G^{p}(\mathbb{D})$. For $\p+\q<0$, we show if  $\partial_z u$ or $\partial_{\overline{z}}u \in H_G^1(\mathbb{D})$ then $u=0$ or $u$ is a polyharmonic function.
\end{abstract}
\section{Introduction }
Let $\D$ be the open unit disc in the complex plane $\C$, let $\T = \partial\D$ be the unit circle, and denote by
$$\partial_z=\frac{1}{2}(\partial_x-i \partial_y)\quad\quad \text{and}\quad\quad \partial_{\overline{z}}=\frac{1}{2}(\partial_x+i \partial_y).$$
\subsection{$(\p,\q)$-harmonic functions}
\noindent Consider the family of second order partial differential operators on $\D$ of the form
\begin{equation*}
\Delta_{\p,\q}:=(1-|z|^2) \left((1-|z|^2)\partial_{z}\partial_{\overline{z}}+\p z\partial_z+\q\,\overline{z}\partial_{\overline{z}}-\p\q\right),
\end{equation*}

\noindent where $\p,\q$ are real parameters. These operators have  been introduced and studied  in higher dimensions by Geller \cite{Ge} and  Ahern et al. \cite{Ah,Ah2} and recently investigated in the planar case  in \cite{KO}.

Of particular interest are the solutions of the associated homogeneous equation in $\D$:
\begin{equation}
\Delta_{\p,\q}\, u=0. \label{Lpq}
\end{equation}

We say that a function $u$ is $(\p,\q)$-harmonic if $u$ is twice continuously differentiable in $\D$ and  $\Delta_{\p,\q} u=0$. We remark if $u$ is $(\p,\q)$-harmonic, then $\overline{u}$ is $(\q,\p)$-harmonic.\\  $(0,\p)$-harmonic functions are called $\alpha$-harmonic functions. \\

An interesting example of an $(\p, \q)$-harmonic function is $K_{\p,\q}$ defined by 
\begin{equation}
    K_{\p,\q}(z)=\frac{(1-|z|^2)^{\p+\q+1}}{(1-z)^{\p+1}(1-\overline{z})^{\q+1}}.\nonumber
\end{equation}

 For $\p,\q \in \R\backslash \Z^-$  such that $\p+\q>-1$, the $(\p,\q)-$harmonic Poisson kernel is defined by 
 \be
P_{\p,\q}(z)=c_{\p,\q}\,\frac{(1-|z|^2)^{\p+\q+1}}{(1-z)^{\p+1}(1-\overline{z})^{\q+1}}, \quad c_{\p,\q}=\frac{\Gamma(\p+1)\Gamma(\q+1)}{\Gamma(\p+\q+1)}, \label{cpq}
\ee

where $\Gamma$ is the Gamma function.\\

The $(\p,\q)-$harmonic Poisson integral of  $f\in L^1(\T)$ is defined by 
$$
P_{\p,\q}[f](z):=\frac{1}{2\pi}\int_{0}^{2\pi} P_{\p,\q}(ze^{-i\theta})f(e^{i\theta})d\theta, \quad z\in\D.$$

 Let 
 $$M_{\p,\q}(z):=P_{\p,\q}[1](z)=\frac{1}{2\pi}\int_{0}^{2\pi} P_{\p,\q}(ze^{-i\theta})\, d\theta.$$
 Then $M_{\p,\q}$ is a radial function and 
 $$
 M_{\p,\q}(z)=c_{\p,\q}F(-\p,-\q;1;|z|^2),
 $$
where $F$ is the Gaussian hypergeometric function. The constant $c_{\p,\q}$ is chosen such that  $\lim_{|z|\to 1} M_{\p,\q}(z)=1$. We should point out that in the case $\p=0$ or $\q=0$, the function $M_{\p,\q}$ is constant and  $M_{\p,\q}(z)=1$. For more properties and H\"older continuity of $(\p,\q)$-harmonic functions, see   \cite{Ah,Xchen, AM, KO, Pei}.\\

Let us recall the Gaussian hypergeometric function defined by
$$
F(a,b;c;z)=\sum_{n\geq 0}\frac{(a)_n(b)_n}{(c)_n}\frac{z^n}{n!},\ \ \ \  z\in\D,
$$
for $a,b,c\in\C$ such that $c\neq -1,-2,...$ where
$$
(a)_n=\begin{cases}
1 &\mbox{for} \ \ n=0,\\
a(a+1)\ldots (a+n-1) &\mbox{for} \ \  n=1,2,3,...

\end{cases}
$$
which are called Pochhammer symbols.\\

We list few properties, see for instance (\cite{And}, Chapter 2)
\begin{eqnarray}
\lim\limits_{x\rightarrow1}F(a,b;c;x)&=&\frac{\Gamma(c)\Gamma(c-a-b)}{\Gamma(c-a)\Gamma(c-b)},\ \ \ \ \ \text{if}\ \ \ c-a-b>0,\label{limitF}\\
F(a,b;c;x)&=&(1-x)^{c-a-b}F(c-a,c-b;c;x),\label{transF}\\
\frac{d}{d x}F(a,b;c;x)&=&\frac{ab}{c}F(1+a,1+b;1+c;x),\label{dF}\\
F(a,b;a+b;x)&\approx& -\frac{\Gamma(a+b)}{\Gamma(a)\Gamma(b)}\log(1-x)  \quad (x\to 1^-)\label{Gauss}.
\end{eqnarray}

 The asymptotic relation (\ref{Gauss}) is due to Gauss, and its refined form is due to  Ramanujan \cite[p. 71]{Be}. \\

Also, we  consider the associated Dirichlet boundary value problem of functions $u$, satisfying the equation (\ref{Lpq}),
\begin{equation}\label{Lpqf}
\begin{cases}
\Delta_{\p,\q}\, u&=0 \quad\text{in}\quad\D,\\
 u&=f\quad\text{on}\quad \T.
\end{cases}
\end{equation}

The  function $f$ is a continuous on $\T$, i.e., $f\in \mathcal{C}(\T)$, and the  condition  (\ref{Lpqf}) is  $\ds\lim_{r\to 1}u_r= f$ in $\mathcal{C}(\T)$,  where $u_r(e^{i\theta}):=u(re^{i\theta}).$\\

\noindent\begin{Thm}[\cite{Ge}]
Let $\p,\q \in \R\backslash \Z^{-}$ be such that $\p+\q > -1$. Let  $f\in\mathcal{C}(\T)$. Then a function $u$ in $\D$ satisfies \eqref{Lpqf} if and only if it has the form
\begin{equation*}
    u(z) =  P_{\p,\q}[f](z)= \frac{1}{2\pi}\int_0^{2\pi}  P_{\p,\q}(ze^{-i\gamma})f(e^{i\gamma})d\gamma,  
\end{equation*}
for any $z\in\D$.
\end{Thm}

The homogeneous series expansion of $(\p,\q)$-harmonic functions is giving by the following theorem.\\

\noindent\begin{Thm}[\cite{Ah}, \cite{KO}]
Let $\p, \q\in\R$. Then $u$ is an $(\p, \q)$-harmonic function if and only if it has the form
\be\label{series expan}
u(z)=\sum_{k\geq 0}c_{k}F(-\p,k-\q;k+1;|z|^2)z^k+\sum_{k\geq 1}c_{-k}F(-\q,k-\p;k+1;|z|^2)\overline{z}^k,
\ee
where $c_k \in \mathbb{C}$, for all $k\in \mathbb{Z}$.
\end{Thm}

Thus for $\p,\q \in \R\backslash \Z^{-}$ with $\p+\q > -1$,   $f\in\mathcal{C}(\T)$, and $u=P_{\p,\q}[f]$, the 
relationship between the coefficients  $c_k$ in \eqref{series expan} and 
$\hat{f}(k)$, the Fourier coefficients of $f$, is given by

\be\label{Fourierck}
c_k F(-\p,k-\q;k+1;1)=\hat{f}(k) \mbox{ 
 for } k\geq 0,
\ee
and 

\be\label{Fourierck2}
c_{-k} F(-\q,k-\p;k+1;1)=\hat{f}(-k) \mbox{ 
 for } k\geq 1.
\ee
\vskip 1em

\noindent Let  $\lambda>-1$ and $r\in[0,1)$, we denote by

$$
I_{\lambda}(r):=\frac{1}{2\pi}\int_{0}^{2\pi} \frac{(1-r^2)^{\lambda}}{|1-re^{-i\theta}|^{\lambda+1}}\, d\theta,
$$

For $\lambda>0$, we have  
\be
\lim_{r\to 1}I_{\lambda}(r)=\frac{\Gamma(\lambda)}{\Gamma^2(\frac{\lambda}{2}+\frac{1}{2})}=:c_{\lambda}.\label{clambda}
\ee

\noindent The following lemma provides some estimates regarding  $I_\lambda$. 

\begin{lem}{\rm (\cite[Proposition 1.1]{KMM} and \cite{oA})}\label{lem1:1}
	Let $r\in[0,1)$.
	\begin{enumerate}
		\item If $\lambda >0$, then 	
		\begin{equation*}
			I_{\lambda} (r) \leq c_{\lambda}.
		\end{equation*}
		\item If $\ -1<\lambda<0$, then
		\begin{equation*}
			I_{\lambda} (r) \leq c_{-\lambda}
			(1-r^2)^{\lambda}.
		\end{equation*}
	\end{enumerate}
	
\end{lem}

 \subsection{Hardy type spaces}
For $p \in(0,\infty]$, {\it the generalized Hardy space} $H_G^p(\D)$ consists of all measurable functions $u$ from $\D$ to $\C$ such that $M_p(r,u)$ exists for all $r \in(0,1)$, and $\|u\|_p < \infty$, where
$$M_p(r,u)=\bigg(\frac{1}{2\pi}\int_0^{2\pi}|u(re^{i\theta})|^p d\theta\bigg)^{\frac{1}{p}}$$
and
$$
\|u\|_p=\begin{cases}
\sup\{M_p(r,u):0<r<1\} &\mbox{if} \ \ p\in(0,\infty) \\
\sup\{|u(z)|:z\in\D\} &\mbox{if} \ \ p=\infty.
\end{cases}$$
 Obviously, $H^p_G(\D) \subset L^p(\D)$ for
each $p\in(0,\infty]$. We denote by $h_{\p,\q}^p(\D)$ the corresponding Hardy space for $(\p,\q)$-harmonic functions.
The classical {\it Hardy space $H^{p}(\mathbb{D})$}  (resp. $h^{p}(\mathbb{D})$) is the set of  all  elements of $H_G^{p}(\mathbb{D})$ which are analytic (resp. harmonic) on $\D$. (cf. \cite{Du,Du1}).
\\

The Hardy theory of $(\p,\q)$-harmonic functions is studied in detail in \cite{Ah}. The authors  proved the following 

\begin{Thm}\cite{Ah}
    Let $u$ be an $(\p,\q)$-harmonic function, and assume $\p+\q>-1$. Then :
    \begin{enumerate}
        \item[(i)] $u=P_{\p,\q}[f]$ for some $f\in L^p(\T)$, $1<p\leq \infty$ if and only if  $\|u\|_p<+\infty$.
        \item[(ii)] $u=P_{\p,\q}[\mu]$ for some measure $\mu$ if and only if  $\|u\|_1<+\infty$.
        
    \end{enumerate}
\end{Thm}

\subsection{Lebesgue spaces} Denote by $L^p(\T)$ where $p\in [1,\infty]$ the space of all measurable functions $f$ of $\T$ into $\C$ with
$$
\|f\|_{L^p}=\begin{cases}
\bigg(\displaystyle \frac{1}{2\pi}\int_0^{2\pi}|f(e^{i\theta})|^p d\theta\bigg)^{\frac{1}{p}} &\mbox{if} \ \ p\in[1,\infty) \\
\text{ess}\sup\{|f(e^{i\theta})|:\theta\in[0,2\pi)\} &\mbox{if} \ \ p=\infty.
\end{cases}$$
It is well-known that if $f$ is absolutely continuous, then it is of bounded variation. This implies that for almost all $e^{i\theta}\in\T$, the derivative $\dot{f}(e^{i\theta})$ exists and integrable on $(0,2\pi)$, where
$$\dot{f}(e^{i\theta}):=\frac{d }{d\theta}f(e^{i\theta}).$$

\subsection{Harmonic conjugates and M. Riesz theorem}
If $u$ is a harmonic  function  on $\D$  given by 
$$u(re^{i\theta})= \sum_{n=-\infty}^{\infty} c_n r^{|n|} e^{in\theta},
$$
then its harmonic conjugate is defined by 
$$\tilde{u}(re^{i\theta}):= \sum_{n=-\infty}^{\infty} m_n c_n r^{|n|} e^{in\theta},$$ 
where $m_n=-i\, {\rm sign}\, n$; in particular $\tilde{u}(0)=0$. In addition, we have  
\be \label{cong}
u(z)+i \tilde{u}(z)=-u(0) + 2 \sum_{n=0}^{\infty} c_n z^n.
\ee
The expression $\sum_{n\geq 0} c_n z^n$ is called the {\it Riesz projection} of $u$, and denoted by $P_+u$ i.e.,
$$P_+ : h(\D) \to H(\D) $$ is  defined by $$P_+(u)(z):=\sum_{n\geq 0} c_n z^n.$$
Thus Eq. (\ref{cong}) can be rewritten as 
\be
u+i\tilde{u}=-u(0)+2P_{+}(u). \nonumber
\ee

It is natural to consider the following question: if $u$ is harmonic in $\D$ and $
 u \in h^p(\D)$, $1\le p \le \infty$, does $P_+(u) \in H^p(\D)$?\\
The question has an affirmative answer for all $1<p<\infty$. This is the content of the famous theorem proved by M. Riesz \cite{Riesz}, see also Rudin \cite[Theorem 17.26]{rud1}.\\

\begin{Thm}(M. Riesz)  For $1<p<\infty$,  the Riesz projection $P_+$ 
$$P_+: h^p(\D) \to H^p(\D) $$
is a bounded operator, i.e.,
there is a constant $A_p$ such that 
\be
\|P_+(u)\|_p \le A_p \| u \|_{p}\nonumber
\ee
for any $\ f \in h^p(\D)$.
\end{Thm}
It turns out that the Riesz theorem  is not true  for  $ p = 1$ and $p = \infty$.\\

In this paper, we consider  the following :\\

\noindent {\bf Question:} Under what conditions on the boundary function $f$ ensure that the partial derivatives of its $(\p,\q)$-harmonic extension $u =P_{\p,\q}[f]$, i.e., $\partial_zu$ and $\partial_{\overline{z}}u$, are in the space $H^p_G(\D)$ (or $L^p(\D)$), where $p\in[1,\infty]$?\\

In \cite{zhu}, Zhu treated the corresponding  problem for  harmonic functions and was improved in \cite{Pou}.

\begin{Thm}[\cite{zhu} and \cite{Pou}]
Suppose that $u = P[f]$ is a harmonic
mapping in $\D$ and $\dot{f} \in L^p(\T)$, where $f$ is an absolutely continuous function.
\begin{itemize}
    \item If $p\in[1,\infty)$, then both $\partial_z u$ and  $\partial_{\overline{z}}u$ are in $L^p(\D)$.
     \item If $p = \infty$, then there exists a harmonic mapping $u = P[f]$, with $\dot{f} \in L^\infty(\T)$, such that neither $\partial_z u$ nor $\partial_{\overline{z}}u$ is in $L^\infty(\D)$.
\end{itemize}
\end{Thm}
Furthermore, by requiring the harmonic extension to be quasiregular \cite{zhu} or more generally elliptic \cite{Pou}, they proved that the partial derivatives are in $H^p(\D)$, for $p \in[1,\infty]$.\\

\begin{Thm}[\cite{Pou}]
Suppose that  $p \in[1,\infty]$ and $u = P[f]$
is a $(K,K')$-elliptic mapping in $\D$ with $\dot{f} \in L^p(\T)$, where $f$ is an absolutely continuous function, $K \geq 1$ and $K'\geq 0$. Then both $\partial_z u$ and $\overline{\partial_{\overline{z}}u}$ are in $H^p(\D)$.
\end{Thm}

In \cite{KM2023-2}, we provide  a refinement of the two previous theorems, we prove that for $1<p<\infty$, both
 $\partial_z u$ and $\overline{\partial_{\overline{z}}u}$ are in $H^{p}(\mathbb{D})$ without any extra conditions on 
$u$.\\

\begin{Thm}
Suppose that $f$ is an  absolutely continuous  function on $\T$ and  $u=P[f]$ is a harmonic mapping in $\mathbb{D}$ and
$\dot{f}\in L^{p}(\mathbb{T})$.
\begin{enumerate}
 
\item
If $p\in(1,\infty)$, then both  $\partial_z u$ and $\overline{\partial_{\overline{z}}u}$ are in $H^{p}(\mathbb{D}).$ Moreover, there exists a constant $A_p$ such that
$$ \max(\|\partial_z u \|_p, \|\partial_{\overline{z}}u\|_p) \leq A_p \| \dot{f}\|_{L^p}.$$

\item If $p=1$ or $p=\infty$, then both  
  $\partial_z u$ and $\overline{\partial_{\overline{z}}u}$ are in $H^p(\mathbb{D})$  if and only if  $H(\dot{f}) \in L^p(\T)$. Moreover
  $$2iz \partial_z u(z)=P[\dot{f}+iH(\dot{f})](z),$$

where $H(\dot{f})$ denotes the Hilbert transform of $\dot{f}$.
\end{enumerate}
\end{Thm}

In \cite{KM2023}, the authors considered the case of $\alpha$-harmonic functions and they proved the following.\\

\begin{Thm}
Let $\alpha \in (-1,\infty)$ with $\alpha\not=0$ and let $u=P_\alpha[f]$ and $f$ is absolutely continuous  such that $\dot{f}\in L^p(\T)$ with $1\leq p\leq \infty$.
\begin{enumerate}
    \item  If $\alpha>0$, then
 $\partial_z u$ and  $\partial_{\overline{z}}u$ are in  $H_G^p(\D)\subset L^p(\D)$.
 \item 
If $\alpha\in (-1,0)$, then  $\partial_z u$ and  $\partial_{\overline{z}}u$ are in $L^p(\D)$ for $p<- \frac{1}{\alpha}$.

\item For $\alpha \in (-1,0)$ and  $p\geq -\frac{1}{\alpha}$ there exists $u$ an  $\alpha$-harmonic function such that $\partial_z u$   and $\partial_{\overline{z}} u$ $\not\in L^p(\D)$, moreover,  $\partial_z u$   and $\partial_{\overline{z}} u$ $\not\in H_G^1(\D)$.

\item Let $\alpha \in (-1,0)$. Then    $\partial_z u$ or $\partial_{\overline{z}}u$  is in  $H_G^p(\D)$, if and only if,   $u$ is analytic.
\end{enumerate}
  \end{Thm}
Recall that $P_\alpha$ is the Poisson kernel for $\alpha$-harmonic functions defined by  $$P_{\p}(z):= P_{0,\p}(z)=\frac{(1-|z|^{2})^{\alpha+1}}{(1-z)(1-\overline{z})^{\alpha+1}}.$$

The aim of this paper is extend the previous results for $(\alpha,\beta)$-harmonic functions for $\p,\q$ real numbers such that $\p,\q\in \R \setminus \Z^-$ with $\p+\q>-1$.

\section{Main results}
Let us recall the polar partial derivative of $u$ with respect to $\theta$ is given by
\be
\partial_\theta u=i(z\partial_zu-\overline{z}\partial_{\overline{z}}u).\label{dtheta1}
\ee
We start this section by the estimate of the angular derivative for $(\p,\q)$-harmonic functions.
\begin{thm}\label{Thm2.1}
Assume that $\p+\q>-1$ and  $p\in [1,\infty]$. Let $u=P_{\p,\q}[f]$ be an $(\p,\q)$-harmonic mapping on $\D$ with the  function $f$ being absolutely continuous with $\dot{f}\in L^p(\T)$. Then 
$$
\partial_\theta u =P_{\p,\q}[\dot{f}]\in h_{\p,\q}^p(\D),
$$ and
\be 
\|\partial_\theta u\|_p\leq |c_{\alpha,\beta}| \, c_{\p+\q+1} \, \| \dot{f}\|_{L^p}, \label{estim dtheta}
\ee

\noindent where the constants $c_{\p,\q}$ and $c_{\p+\q+1}$ are defined in  (\ref{cpq}) and (\ref{clambda}) respectively.

\end{thm}
\begin{pf} By assumptions, we have

$$u(z) = P_{\p,\q}[f](z)= \frac{c_{\p,\q}}{2\pi}\int_0^{2\pi} K_{\p,\q}(ze^{-it} )\, f(e^{it})\, dt.$$
Let $z=r e^{i\theta}$, an easy computations show that 
 $$\partial_{\theta} K_{\p,\q} (ze^{-it})= -\partial_{t} K_{\p,\q} (ze^{-it}).$$

\noindent An integration by parts leads to 
\be
\partial_\theta u(z)=P_{\p,\q}[\dot{f}](z),\nonumber
\ee
\noindent By  Jensen's inequality, we get
$$|\partial_\theta u(z)|^p\leq (|c_{\p,\q}|\, I_{\p+\q+1}(r))^{p-1}\, \frac{1}{2\pi}\int_0^{2\pi}|K_{\p,\q}(ze^{-it})||\dot{f}(e^{it})|^p\, dt.$$
Using Fubini's theorem, it yields
$$\frac{1}{2\pi}\int_{0}^{2\pi}|\partial_\theta u(re^{i\theta})|^p d\theta\leq (|c_{\p,\q}|\, I_{\p+\q+1}(r))^{p}\| \dot{f}\|_{L^p}^p. $$
Thus, using Lemma \ref{lem1:1}, we have 
$$\|\partial_\theta  u\|_p \leq |c_{\p,\q}|\, c_{\p+\q+1}  \| \dot{f}\|_{L^p}.$$
For $p=\infty$, it follows that
$$|\partial_\theta u(z)|\leq |c_{\p,\q}| \, |I_{\p+\q+1}(r)\|\dot{f}\|_{L^\infty}\leq |c_{\p,\q}|c_{\p+\q+1 }\| \dot{f}\|_{L^\infty}.$$

\end{pf}

\begin{lem}
\be\label{dKpq}
\partial_z K_{\p,\q}(z)=\frac{(\p+1)(1-|z|^2)^{\p+\q}}{(1-z)^{\p+2}(1-\overline{z})^{\q}}-\frac{\q\,\overline{z}(1-|z|^2)^{\p+\q}}{(1-z)^{\p+1}(1-\overline{z})^{\q+1}},
\ee
and
\be\label{dbarKpq}
\partial_{\overline{z}} K_{\p,\q}(z)=\frac{(\q+1)(1-|z|^2)^{\p+\q}}{(1-z)^{\p}(1-\overline{z})^{\q+2}}-\frac{\p\,z(1-|z|^2)^{\p+\q}}{(1-z)^{\p+1}(1-\overline{z})^{\q+1}}.
\ee
\end{lem}

\begin{pf}
  Differentiating the kernel $K_{\p,q}$ with respect to $z$, we get
  
  	$$\partial_z K_{\p,\q}(z)=\bigg(-(\p+\q+1)\frac{\overline{z}}{1-|z|^2}+(\p+1)\frac{1}{1-z}\bigg)K_{\p,\q}(z),$$
  see for instance \cite[Lemma 1.1]{KO}. Next, observe that 
   $$-(\p+\q+1)\frac{\overline{z}}{1-|z|^2}+(\p+1)\frac{1}{1-z}=\frac{(\p+1)(1-\overline{z})}{(1-z)(1-|z|^2)}-\frac{\q \overline{z}}{(1-|z|^2)},
   $$
and the conclusion \eqref{dKpq} follows immediately. By the symmetry property of the kernel $K_{\p,\q}$, we get \eqref{dbarKpq}. 
\end{pf}

\begin{prop}
Let  $u=P_{\p,\q}[f]$ be an  $(\p,\q)$-harmonic mapping on $\D$ with $f$ being absolutely continuous on $\T$. Then
\begin{equation}
z\partial_z u= -ic_{\p,\q} \left(K_{\p.\q-1}[\dot{f}]+i\q\overline{z}\, K_{\p-1,\q}[f_1]\right),\label{zdzu}
\end{equation}
\be
\overline{z}\partial_{\overline z} u= i c_{\p,\q} \left( K_{\p-1.\q}[\dot{f}]-i\p zK_{\p,\q-1}[f_1]\right),\label{zbaru}
\ee

where $f_1$ is defined on $\T$ by $f_1(e^{it})=e^{it} f(e^{it})$.
\end{prop}

\begin{pf}
The partial derivative of $u$ with respect to $z$ is giving by 
$$\partial_z u(z)=\frac{1}{2\pi}\int_0^{2\pi}\partial_z P_{\p,\q}(ze^{-it}) e^{-it} f(e^{it})dt.$$
By \eqref{dKpq}, we have 
$$\partial_z K_{\p,\q}(z)=\frac{(\p+1)(1-|z|^2)^{\p+\q}}{(1-z)^{\p+2}(1-\overline{z})^{\q}}-\frac{\q\,\overline{z}(1-|z|^2)^{\p+\q}}{(1-z)^{\p+1}(1-\overline{z})^{\q+1}}.$$
Hence 

\begin{eqnarray}
\frac{1}{c_{\p,\q}}\partial_z u(z)&=&\frac{1}{2\pi}(1-|z|^2)^{\p+\q}\int_0^{2\pi}\frac{(\p+1)f(e^{it})e^{-it}}{(1-ze^{-it})^{\p+2}(1-\overline{z}e^{it})^{\q}}\, dt\nonumber\\
&-&\frac{1}{2\pi}(1-|z|^2)^{\p+\q}\int_0^{2\pi}\frac{\q\overline{z}f(e^{it})}{(1-ze^{-it})^{\p+1}(1-\overline{z}e^{it})^{\q+1}}\, dt.\nonumber
\end{eqnarray}

\noindent Introduce the functions $h_1$ and $h_2$ defined on $\T$ by  
$$h_1(e^{it}):=\frac{(\p+1)e^{-it}}{(1-ze^{-it})^{\p+2}(1-\overline{z}e^{it})^{\q}}$$
and
$$h_2(e^{it}):=\frac{\q\overline{z}}{(1-ze^{-it})^{\p+1}(1-\overline{z}e^{it})^{\q+1}},$$
so that we get the following expression
\be
\frac{1}{c_{\p,\q}} \partial_z u(z)=\frac{1}{2\pi}(1-|z|^2)^{\p+\q}\bigg(\int_0^{2\pi}h_1(e^{it})f(e^{it})dt -\int_0^{2\pi}h_2(e^{it})f(e^{it})dt \bigg).\label{dzu}
\ee
Consider the function $g$ on $\T$ defined by  
$$g(e^{it}):=\frac{1}{(1-ze^{-it})^{\p+1}(1-\overline{z}e^{it})^{\q}},$$

\noindent its derivative $\dot{g}$ is given by 
$$\dot{g}(e^{it})=\frac{-ize^{-it}(\p+1)}{(1-ze^{-it})^{\p+2}(1-\overline{z}e^{it})^{\q}}+\frac{i\overline{z}\q e^{it}}{(1-ze^{-it})^{\p+1}(1-\overline{z}e^{it})^{\q+1}},$$
that is 
\be 
\dot{g}(e^{it}) \frac{}{}=-izh_1(e^{it})+ie^{it}h_2(e^{it}).\label{gdot}
\ee

\noindent Using \eqref{dzu} and \eqref{gdot}, we have

\begin{eqnarray}
\frac{iz\partial_z u(z)}{c_{\p,\q} (1-|z|^2)^{\p+\q}}&=&\frac{1}{2\pi}\int_0^{2\pi}\big(-\dot{g}(e^{it})+ie^{it}h_2(e^{it})\big)f(e^{it})\,dt -\frac{1}{2\pi}\int_0^{2\pi}iz h_2(e^{it})f(e^{it})\, dt\nonumber\\
&=&-\frac{1}{2\pi}\int_0^{2\pi}\dot{g}(e^{it})f(e^{it})dt +\frac{i}{2\pi}\int_0^{2\pi}\big(e^{it}-z\big)h_2(e^{it})f(e^{it})dt\nonumber
\end{eqnarray}
By integration by parts in the first integral and using the fact that $e^{it}-z=e^{it}(1-ze^{-it})$, we obtain
\begin{eqnarray}
\frac{iz\partial_z u(z)}{c_{\p,\q}(1-|z|^2)^{\p+\q}}&=&\frac{1}{2\pi}\int_0^{2\pi}\frac{\dot{f}(e^{it})}{(1-ze^{-it})^{\p+1}(1-\overline{z}e^{it})^{\q}}\,dt\nonumber\\ &+&\frac{i}{2\pi}\int_0^{2\pi}\frac{\q\overline{z}e^{it}f(e^{it})}{(1-ze^{-it})^{\p}(1-\overline{z}e^{it})^{\q+1}}\, dt.\nonumber
\end{eqnarray}
Finally, we deduce that

$$iz\partial_z u= c_{\p,\q} \left( K_{\p.\q-1}[\dot{f}]+i\q\overline{z}K_{\p-1,\q}[f_1] \right).$$
The equation \eqref{zbaru} is obtained by considering $\overline{u}$ which is $(\q,\p)$-harmonic. \end{pf}

Now, we are ready to state our first main  result.

\begin{thm}\label{2.2}
Assume that $\p+\q>-1$ and $1\leq p \leq \infty$. Let $u=P_{\p,\q}[f]$ be an $(\p,\q)$-harmonic
mapping on $\D$ with $f$ being absolutely continuous and  $\dot{f}\in L^p(\T)$. Then
\be\label{2...8}
M_p(r,z\partial_z u)\leq |c_{\p,\q}|\, I_{\p+\q}(|z|)\bigg(\|\dot{f}\|_{L^p}+|\q|\|f\|_{L^p}\bigg),
\ee
and 
\be\label{2..9}
M_p(r,\overline{z}\partial_{\overline{z}} u)\leq |c_{\p,\q}|\, I_{\p+\q}(|z|)\bigg(\|\dot{f}\|_{L^p}+|\p|\|f\|_{L^p}\bigg).
\ee

\noindent In particular 
\begin{enumerate}
    \item[(i)] if $\p+\q>0$, then $\partial_z u$ and $\partial_{\overline{z}} u$ belong to $H_G^p(\D)$ and 
\be\label{2...10}
\|z\partial_z u\|_p\leq  |c_{\p,\q}|\, c_{\p+\q} \big(\|\dot{f}\|_{L^p}+|\q|\|f\|_{L^p}\big),
\ee
and 
 
\be\label{2...11}
\|\overline{z}\partial_{\overline{z}} u\|_p\leq  |c_{\p,\q}|\, c_{\p+\q} \big(\|\dot{f}\|_{L^p}+|\p|\|f\|_{L^p}\big).
\ee
\item [(ii)] If $\p+\q <0$ and $p<-\frac{1}{\p+\q}$, then $\partial_z u$ and $\partial_{\overline{z}} u$ are in $L^p(\D)$.

\item[(iii)] If $\p+\q <0$ and  $p\geq -\frac{1}{\p+\q}$, then there exists $u$ an  $(\p,\q)$-harmonic function such that $\partial_z u$   and $\partial_{\overline{z}} u$ $\not\in L^p(\D)$. Moreover,  $\partial_z u$   and $\partial_{\overline{z}} u \not\in \mathcal{H}_{\mathcal{G}}^1(\D)$.
\end{enumerate}
\end{thm}
\begin{pf}
Using \eqref{zdzu}, we obtain 
$$
M_p(r,z\partial_z u)\leq |c_{\p,\q}|  \left( M_p(r,K_{\p.\q-1}[\dot{f}]) +|\q||z|M_p(r,K_{\p-1,\q}[f_1]) \right).
$$
Similarly, as in the proof of Theorem \ref{Thm2.1}, using Jensen's inequality and Fubini's theorem, we have
$$
M_p(r,K_{\p.\q-1}[\dot{f}]) \leq I_{\p+\q}(|z|) \|\dot{f}\|_{L^p} \quad \mbox{and} \quad M_p(r,K_{\p-1.\q}[f_1]) \leq I_{\p+\q}(|z|) \|f\|_{L^p}.
$$
Therefore, we obtain the desired result \eqref{2...8}.  The equations \eqref{2...10}, \eqref{2...11} and the assertion (ii) are immediate consequences of \eqref{2...8},  \eqref{2..9} and Lemma  \ref{lem1:1}.\\

(iii) Let $k\in \Z$ and 
 $$
 M_{\p,\q,k}(z):= P_{\p,\q}[e^{ik\theta}](z)=
 \frac{1}{2\pi}\int_{0}^{2\pi} P_{\p,\q}(ze^{-i\theta})e^{i k\theta}d\theta, \quad z\in\D.
 $$
From the identity, see for instance \cite{KO}
 $$\begin{aligned} \frac{(1-|z|^2)^{\p+\q+1}}{(1-z)^{\p+1} (1-{\bar{z}})^{\q+1}}&=\sum _{m=0}^\infty \frac{(\p+1)_m}{m!}F(-\p,m-\q;m+1;|z|^2)z^m\\&\quad \ +\sum _{m=1}^\infty \frac{(\q+1)_m}{m!}F(-\q,m-\p;m+1;|z|^2){\bar{z}}^m, \end{aligned}$$
 
 it yields for $ k\in \mathbb{N}$
 \be 
 M_{\p,\q,k}(z)= c_{\p,\q}\frac{(\p+1)_k}{k!}\, F(-\p,k-\q;k+1;|z|^2)z^k.\nonumber
 \ee
Differentiate $M_{\p,\q,k}$ with respect to $\overline{z}$,  using (\ref{dF}) and \eqref{transF}, we obtain 
\begin{eqnarray}
\partial_{\bar{z}} M_{\p,\q,k}(z)&=&- c_{\p,\q}\frac{(k-\q)(\p)_{k+1}}{(k+1)!}\, F(1-\p,1+k-\q;k+2;|z|^2)z^{k+1}\nonumber \\
 &=&- c_{\p,\q}\frac{(k-\q)(\p)_{k+1}}{(k+1)!}\, (1-|z|^2)^{\p+\q}F(k+\p+1,\q+1;k+2;|z|^2)z^{k+1}.\nonumber
 \end{eqnarray}
By (\ref{limitF}), remark that 
$$\lim_{|z|\to 1} F(k+\p+1,1+\q;k+2;|z|^2)= \frac{\Gamma(k+2) \Gamma(-\p-\q)}{\Gamma(1-\p) \Gamma(k+1-\q)}.
$$

Next, we will pick $k\not = \q$ to ensure that $\lim_{|z|\to 1} F(k+\p+1,1+\q;k+2;|z|^2)$ is {\it non-zero}, that is, $1-\p$ and $k+1-\q$ are not negative integers. In that case, we have
\be
|\partial_{\bar{z}} M_{\p,\q,k}(z)| \approx  (1-|z|^2)^{\p+\q} \quad (|z|\to 1^-).\nonumber
\ee
Thus $\partial_{\bar{z}} M_{\p,\q,k} \not\in L^p(\D)$ for $p\geq -\frac{1}{\p+\q}$.  In addition, since $\p+\q<0$, we have  $\partial_{\bar{z}} M_{\p,\q,k} \not \in H^1_G(\D) $.\\

 For the choice of $k$, we will consider the following three cases:
 
 \begin{enumerate}
      \item[(i)] For $\q\in \mathbb{N}$,  take $k=\q+1$. We have $1-\p>1+\q\geq 1$, we see that $\lim_{|z|\to 1} F(\p+\q+1,\q;\q+2;|z|^2)$ is non-zero. 
 \item[(ii)] The case  $\p\in \mathbb{N}$ is treated similarly.
 \item[(iii)] For both $\p$ and $\q$ are non integers, take $k=0$. Trivially, as $\q\not \in \mathbb{N}$, we have  $1-\q\not \in \Z^-$. 
 \end{enumerate}
\end{pf}

These examples are introduced in \cite{AM} to  show to failure of the stability property of being  Lipschitz continuous for $(\p,\q)$-harmonic functions when $\p+\q \in (-1,0)$.\\

\subsection{Case $-1<\p+\q<0$}\hfill\\

\noindent  In this section, we  give necessary  conditions  on $u$ an $(\p,\q)$-harmonic  such that   $\partial_z u$ or $\partial_{\overline{z}} u$ is in $H_G^1(\D)$. We show that $u$ is zero or $u$ is  polyharmonic.\\

\noindent As $u$ is $(\p,\q)$-harmonic,  we consider its series expansion, see \eqref{series expan}
$$u(z)=\sum_{k\geq 0}c_{k}F(-\p,k-\q;k+1;|z|^2)z^k+\sum_{k\geq 1}c_{-k}F(-\q,k-\p;k+1;|z|^2)\overline{z}^k.$$
For $n\in  \Z$, let 
\be
u_n(z):=\frac{1}{2\pi}\int_0^{2\pi}u(ze^{i\theta})e^{-in\theta}d\theta.\label{2.8}
\ee

\noindent Thus
\be\label{2.9}
u_{-n}(z)=c_{-n}F(-\q,n-\p;n+1;|z|^2)\overline{z}^n, \quad \mbox{ for } n\geq 1,
\ee
and  
\be\label{2.10}
u_{n}(z)=c_{n}F(-\p,n-\q;n+1;|z|^2)z^n, \quad  \mbox{ for } n \geq 0.
\ee

\noindent By differentiation under the integral sign in (\ref{2.8}),  we obtain
\be
\partial_z u_n(z)=\frac{1}{2\pi}\int_0^{2\pi}\partial_zu(ze^{i\theta})e^{i(1-n)\theta}d\theta,\label{dzuk3}
\ee    

\be
\partial_{\overline z} u_n(z)=\frac{1}{2\pi}\int_0^{2\pi}\partial_{\overline z}u(ze^{i\theta})e^{-i(1+n)\theta}d\theta.\nonumber
\ee

 \noindent Differentiating  the equations (\ref{2.9}) and (\ref{2.10}) with respect to $z$, and using \eqref{transF}, we get  

\be\label{2.13}
\partial_z u_{-n}(z)=c_{-n}\frac{(-\q)(n-\p)}{n+1}(1-|z|^2)^{\p+\q}F(n+\q+1,\p+1;n+2;|z|^2)\overline{z}^{n+1},
\ee 

\begin{eqnarray}\label{2.14}
\partial_z u_{n}(z)&=&c_{n}\frac{(-\p)(n-\q)}{n+1}(1-|z|^2)^{\p+\q}F(n+\p+1,\q+1;n+2;|z|^2)\overline{z}z^n\nonumber\\
&+&
nc_n  F(-\p,n-\q;n+1,|z|^2)z^{n-1}.
\end{eqnarray}

\begin{thm}\label{2.3}
Let $\p$ and $\q$ be non-integers such that $-1<\p+\q<0$ and let $u$ be an $(\p,\q)$-harmonic function.
 If   $\partial_z u$ or $\partial_{\overline{z}}  u$ is in $H_G^1(\D)$, then  $u$ is  zero.
\end{thm}

\begin{pf} Assume that $\partial_z u \in  H_G^1(\D)$. Using (\ref{dzuk3}), we conclude that   $\partial_z u_n$ is {\it bounded} on $\D$ for all $n\in \Z$. \\
On  one hand, according to (\ref{2.13}), for $n\geq 1$,  we have
\be
\partial_z u_{-n}(z)=c_{-n}\frac{(-\q)(n-\p)}{n+1}(1-|z|^2)^{\p+\q}F(n+\q+1,\p+1;n+2;|z|^2)\overline{z}^{n+1}, \nonumber
\ee 
and by the property \eqref{limitF}, we have
\be
\lim_{|z|\to 1}F(n+\q+1,\p+1;n+2;|z|^2)=\frac{\Gamma(n+2)\Gamma(-\p-\q)}{\Gamma(1-\q)\Gamma(n+1-\p)}.\label{2..16}
\ee

As $\p$ and $\q$ are non-integers, we deduce that the limit in (\ref{2..16}) is {\it non-zero}, therefore $\partial_z u_{-n}$ cannot be bounded unless $c_{-n}=0$, for $n\geq 1$.

On the other hand, according to (\ref{2.14}), the partial derivative with respect to $z$ of $u_n$ is given by
\begin{eqnarray}
\partial_z u_{n}(z)&=&c_{n}\frac{(-\p)(n-\q)}{n+1}(1-|z|^2)^{\p+\q}F(n+\p+1,\q+1;n+2;|z|^2)\overline{z}z^n\nonumber\\
&+&
nc_n  F(-\p,n-\q;n+1,|z|^2)z^{n-1}.\nonumber
\end{eqnarray}

As $\p+\q>-1$, remark that the function $F(-\p,n-\q;n+1;|z|^2)$ has a finite limit as $|z|$ goes to $1$, therefore it is bounded on $\D$. In addition notice that   $\ds \lim_{|z|\to 1} F(n+\p+1,\q+1;n+2;|z|^2)$ is {\it non-zero}, as in (\ref{2..16}).
Therefore, since  $\partial_z u_{n}$ is bounded, it yields  $c_n=0$, for $n\geq 0$. Finally, we conclude that $u=0$. \\
The case $\partial_{\overline{z}}  u \in H_G^1(\D)$ is treated similarly and we conclude that $u$ must be zero.
\end{pf}

\begin{thm}
Let $\p$  be  a positive integer and $\q$ be a real number  such that  $-1<\p+\q<0$. Let $u$ be an $(\p,\q)$-harmonic mapping.\\
If  $\partial_z u$ or $\partial_{\overline{z}}  u$ is in $H_G^1(\D)$, then  $c_{-n}=0$ for $n> \p$, where $c_{-n}$ are coefficients in the series expansion of $u$ in  (\ref{series expan}). \\   Indeed, $u$ is   $(\p+1)$-harmonic, and 
$$
u(z)= K_0(z)+K_1(z)|z|^2+\ldots+K_{\p}(z)|z|^{2\p},
$$

where $K_i=h_i +\overline{g_i}$  for $i=0,\ldots,\p-1$ and $K_{\p}$ is analytic, and  $g_i$ are polynomials  and $h_i$ are analytic for $i=0,\ldots, \p-1$.
\end{thm}

Let $n$ be a  positive integer, recall that a function $u$ is called $n$-harmonic if $\Delta^n u=0$.

\begin{pf}

(a) Assume that $\partial_z u\in H^1_G(\D)$, then we deduce from \eqref{dzuk3} that  $\partial_z u_n$ is bounded for all $n\in \Z$. Using (\ref{2.13}), we see that the limit in (\ref{2..16}) cannot be zero for $n>\p$ as by assumptions $\q$ is non zero and $1-\q$ cannot be zero or a negative integer. 
Hence we conclude that 
$c_{-n}=0$ for $n> \p$. Therefore 
$$u(z)=\sum_{k\geq 0}c_{k}F(-\p,k-\q;k+1;|z|^2)z^k+\sum_{k=1}^\p c_{-k}F(-\q,k-\p;k+1;|z|^2)\overline{z}^k.$$

We should notice that the functions $F(-\p,k-\q;k+1;x)$ and $F(-\q,k-\p;k+1;x)$ are polynomials respectively of degree $\p$ and $\p-k$ ($1\leq k\leq \p$), and we have 

\be
\sum_{k\geq 0}c_{k}F(-\p,k-\q;k+1;|z|^2)z^k= h_0(z)+h_1(z)|z|^2+\ldots+h_{\p}(z)|z|^{2\p},
\nonumber
\ee
where $h_i$ are analytic for $i=0,\ldots,\p$. In addition,

\be
\sum_{k=1}^\p c_{-k}F(-\q,k-\p;k+1;|z|^2)\overline{z}^k=\overline{g_0(z)}+ \overline{g_1(z)} |z|^2 +\ldots + \overline{g_{\p-1}(z)}\, |z|^{2(\p-1)}, 
\nonumber
\ee
\noindent where $g_i$ are polynomials  for $i=0,\ldots, \p-1$.  Therefore, we obtain the desired result.  


(b)  Assume that $\partial_{\overline z} u\in H^1_G(\D)$. As in the previous case, we deduce that   $\partial_{\overline z} u_{-n}$ is bounded for $n\geq 1$. 
\begin{eqnarray}
\partial_{\overline z} u_{-n}(z)&=& c_{-n} \frac{\q(\p-n)}{n+1} (1-|z|^2)^{\p+\q}F(n+\q+1,1+\p;n+2;|z|^2)\overline{z}^n z \nonumber\\
&+&nc_{-n} F(-\q,n-\p;n+1;|z|^2)\overline{z}^{n-1}.\nonumber
\end{eqnarray}

As $\p+\q>-1$, remark that the function $F(-\q,n-\p;n+1;|z|^2)$ has a finite limit as $|z|$ goes to $1$, therefore it is bounded on $\D$. 
Similarly as in the previous case 
$\ds \lim_{|z|\to 1} F(n+\q+1,\p+1;n+2;|z|^2)$ is non-zero, for $n> \p$.
Therefore  $c_{-n}=0$, for $n> \p$. 
\end{pf}

For the converse, we prove the following

\begin{thm}\label{converse}
    Let $\p$ be a positive integer and $\q$ be a real number such that $\p+\q \in (-1,0)$.  Assume $f\in \mathcal{C}^{\p-1}(\T)$ with $f^{(\p-1)}$  absolutely continuous  and $f^{(\p)}\in L^p(\T)$ for some $p\in(1,\infty)$ and $ u=P_{\p,\q}[f]$ with $c_{-n}=0$ for $n>\p$.  
    Then $\partial_z u$ and $\partial_{\overline{z}}  u$ are in $H^p_G(\D)$.
\end{thm}

Here $f^{(0)}=f$ and $f^{(\p)}(e^{i\theta}):= \frac{d}{d\theta} f^{(\p-1)}(e^{i\theta})$.
\begin{pf}
   As $c_{-n}=0$ for $n>\p$,  $u$ can be written as  
$$u(z)=\sum_{n\geq 0}c_{n}F(-\p,n-\q;n+1;|z|^2)z^n+\sum_{n=1}^\p c_{-n}F(-\q,n-\p;n+1;|z|^2)\overline{z}^n.$$

\noindent Remark that for $n \in [\![1..\p]\!]$, the hypergeometric function $F(-\q,n-\p;n+1;|z|^2)$ is a polynomial, therefore  the mapping $z\mapsto \sum_{n=1}^\p c_{-n}F(-\q,n-\p;n+1;|z|^2)\overline{z}^n$ is  smooth on $\overline{\D}$. Therefore $\partial_{\overline{z}} u\in H^p(\D)$ if and only  $\partial_{\overline{z}} u_1 \in  H^p(\D)$, where
$$ u_1(z)= \sum_{n\geq 0}c_{n}F(-\p,n-\q;n+1;|z|^2)z^n.$$

Since
$$F(-\p, n-\q;n+1;|z|^2)= \sum_{k=0}^\p  \frac{(-\p)_k (n-\q)_k}{(n+1)_k k!} |z|^{2k},$$

\noindent  we write 
$$u_1(z)= H_0(z)+H_1(z)|z|^2+\ldots+H_{\p}(z)|z|^{2\p},
$$
where for $k\in [\![0..\p]\!]$,  $H_k$ is analytic and  $$H_k(z)= \frac{(-\p)_k}{k!} \sum_{n\geq 0} \frac{(n-\q)_k}{(n+1)_k} c_n z^n. 
$$


\noindent In addition, by \eqref{Fourierck}, we have 

$$c_n\,  \frac{\Gamma(n+1) \Gamma(\p+\q+1)}{\Gamma(\q+1) \Gamma(n+\p+1)}=\hat{f}(n)  \mbox{ for } n\geq 0.
$$
Thus, as 
$$\frac{(n+1)_{\p}}{(n+1)_k}= (n+k+1)_{\p-k} \mbox{ for } k\in[\![0.. \p]\!],$$ we get

 $$H_k(z)= \frac{(-\p)_k \Gamma(\q+1)}{ \Gamma(\p+\q+1)k!} \sum_{n\geq 0} (n-\q)_k (n+1+k)_{\p-k} \, \hat{f}(n) z^n,
$$
and $H_k$ is a linear  combination of $Q_\ell(z):=\ds\sum_{n\geq 0} (in)^\ell \hat{f}(n) z^n$ for $\ell=0,\ldots, \p$. 

We claim that 

$$Q_\ell \in   \bigcap_{1<q<\infty} H^q(\D) \mbox{ for } \ell=0\ldots \p-1 \mbox{ and }  Q_\p \in H^{p}(\D).$$

\noindent Indeed, consider  the harmonic function defined as the Poisson extension of $f$ i.e., $$U(z)=P[f]= \sum_{n\geq 0} \hat{f}(n) z^n +\sum_{n\geq 1} \hat{f}(-n)\overline{z}^n.$$  
\noindent For $\ell=1,\ldots, \p$, we have 
 $$\frac{d^\ell}{ d\theta^\ell} U(z)=P[f^{\ell}]=\sum_{n\geq 0} (in)^\ell\hat{f}(n) z^n +\sum_{n\geq 1} (-in)^\ell\hat{f}(-n)\overline{z}^n.$$

\noindent Hence $Q_{\ell}=P_+(\frac{d^\ell U}  { d\theta^\ell} ),$ the Riesz projection of the harmonic function $\frac{d^\ell U}  { d\theta^\ell} $. Thus the claim follows from Riesz's  theorem on conjugate functions.
Since
    $\partial_{\overline{z}} u_1=zH_1(z)+\ldots + \p z^{\p}H_{\p}(z)\overline{z}^{\p-1},$
this implies that $\partial_{\overline{z}} u \in H^p_G(\D)$.\\ On the other hand, as $\partial_\theta u=i(z\partial_zu-\overline{z}\partial_{\overline{z}}u) \in H^p_G(\D)$, we conclude that $\partial_z u\in H^p_G(\D)$.\end{pf}

We  expect a stronger version of Theorem \ref{converse} by relaxing the assumptions of the boundary function $f$. For $\alpha=1$, it yields  

\begin{cor}
    Let $\beta \in (-2,-1)$ and $p\in (1,\infty)$ and $u=P_{1,\beta}[f]$ where $f$ is absolutely continuous on $\T$ with $\dot{f} \in L^p(\T)$.\\ 
    Then   $\partial_z u$ or $\partial_{\overline{z}}  u$ is in $H_G^1(\D)$, if and only if $\hat{f}(-n)=0$ for $n>1$. In particular $u$ is biharmonic.
\end{cor}
\vskip 1em

\subsection{Case $\p+\q=0$}\hfill\\

\begin{thm}
Let $\p,\q$ be non zero numbers in $\R\setminus \Z^-$ such that $\p+\q=0$ and let $u$ be a $(\p,\q)$-harmonic function.
 If   $\partial_z u$ or $\partial_{\overline{z}}  u$ is in $H_G^1(\D)$, then  $u$ is  zero.
\end{thm}

\begin{pf}
    Assume that $\partial_z u \in  H_G^1(\D)$. Using (\ref{2.8}), we conclude that   $\partial_z u_n$ is {\it bounded} on $\D$. 
According to (\ref{2.13}), we have
\be
\partial_z u_{-n}(z)=c_{-n}\frac{\p(n-\p)}{n+1}F(n-\p+1,\p+1;n+2;|z|^2)\overline{z}^{n+1},
\nonumber
\ee 
and by \eqref{Gauss}, we have
\be
F(n-\p+1,\p+1;n+2;x) \approx \frac{\Gamma(n+2)}{\Gamma(n-\p+1)\Gamma(\p+1)}\log(1-x), \quad (x\to 1^-).\label{2.16}
\ee

As $\p$ is not an integer, we deduce that the limit of the hypergeometric function  in (\ref{2.16}) is  $\infty$ as $x\to 1$, therefore $\partial_z u_{-n}$ cannot be bounded unless $c_{-n}=0$, for $n\geq 1$.

Moreover, according to (\ref{2.14}), the partial derivative with respect to $z$ of $u_n$ is given by
\begin{eqnarray}
\partial_z u_{n}(z)&=&c_{n}\frac{(-\p)(n+\p)}{n+1}F(n+\p+1,1-\p;n+2;|z|^2)\overline{z}z^n\nonumber\\
&+&
nc_n  F(-\p,n+\p;n+1,|z|^2)z^{n-1}.\nonumber
\end{eqnarray}

\noindent Remark that the function $F(-\p,n+\p;n+1;|z|^2)$ has a finite limit as $|z|$ goes to $1$, therefore it is bounded on $\D$. In addition notice that, as in (\ref{2.16}),   $\lim_{|z|\to 1} F(n+\p+1,1-\p;n+2;|z|^2)=\infty$.
Therefore $\partial_z u_{n}$ is bounded unless $c_n=0$, for $n\geq 0$. Finally, we conclude that $u=0$.
\end{pf}

Finally,  given $u=P_{\p,\q}[f]$ a quasi-regular $(\p,\q)$-harmonic function with $\dot{f} \in L^p(\T)$, we prove that  both $\partial_z u$ and  $\partial_{\overline{z}}  u$ are in $H_G^p(\D)$ for $p\in[1,\infty]$.\\

 \noindent Recall that 	a mapping $u:\D \to \C$ is said to be $K-$quasi-regular,
	if $$|\partial_z u(z)|+|\partial_{\overline{z}}u(z)| \leq K (|\partial_z u(z)|-|\partial_{\overline{z}}u(z)|),$$ on $\D$, where $K \in[1, \infty)$ is a constant.
\begin{thm}
    Suppose $1 \leq p \leq  \infty$,  $\p,\q \in \R \setminus \Z^-$ such that  $\p+\q>-1$ and  $u = P_{\p,\q}[f]$ is an $(\p,\q)$-harmonic quasi-regular mapping of
$\D$ with the boundary function $f$ is absolutely continuous and satisfies $\dot{f} \in L^p(\T )$
Then $\partial_z u$ and  $\partial_{\overline{z}}  u$ are in $H_G^p(\D)$.
\end{thm}

\begin{pf}
Let $z=re^{i\theta}$, by using the identity \eqref{dtheta1} and the inequality \eqref{estim dtheta}, we get
$$\frac{1}{2\pi}\int_0^{2\pi}\big(|z\partial_z u(re^{i\theta})|-|\overline{z}\partial_{\overline{z}} u(re^{i\theta})|\big)^p d\theta\leq\frac{1}{2\pi}\int_0^{2\pi}|\partial_{\theta} u(re^{i\theta})|^p d\theta\leq (|c_{\p,\q}|\, c_{\p+\q+1})^p\|\dot{f}\|^p_{L^{p}},$$
where $1\leq p<\infty.$ As  $u$ is quasi-regular, it follows that there exists a constant $K\geq 1$, such that
$$|\partial_z u(z)|+|\partial_{\overline{z}} u(z)|\leq K \big(|\partial_z u(z)|-|\partial_{\overline{z}} u(z)|\big).$$
Hence 
$$\frac{1}{2\pi}\int_0^{2\pi}\big(|z\partial_z u(re^{i\theta})|+|\overline{z}\partial_{\overline{z}} u(re^{i\theta})|\big)^p d\theta\leq (K|c_{\p,\q}|c_{\p+\q+1})^p\|\dot{f}\|^p_{L^{p}},$$
it yields that
$$\frac{1}{2\pi}\int_0^{2\pi}|z\partial_z u(re^{i\theta})|^p d\theta\leq  (K|c_{\p,\q}|c_{\p+\q+1})^p\|\dot{f}\|^p_{L^{p}}$$
and 
$$\frac{1}{2\pi}\int_0^{2\pi}|\overline{z}\partial_{\overline{z}} u(re^{i\theta})|^p d\theta\leq \bigg(\frac{K-1}{2} |c_{\p,\q}|c_{\p+\q+1}\bigg)^p \|\dot{f}\|^p_{L^{p}}.$$
Thus 
$$M_p(r,z\partial_z u)\leq K|c_{\p,\q}|c_{\p+\q+1}\|\dot{f}\|_{L^{p}}\quad\quad\text{and}\quad\quad M_p(r,\overline{z}\partial_{\overline{z}}u)\leq \frac{K-1}{2}c_{\p+\q+1}|c_{\p,\q}|\|\dot{f}\|_{L^{p}}.$$
Hence $\partial_z u$ and  $\partial_{\overline{z}}  u$ are in $H_G^p(\D)$.
\end{pf}

\section{$T_\alpha$-harmonic functions}

 Let $\p \in \R$ and denote by 
 $$\omega_\p(z)=(1-|z|^2)^\p; \quad z\in \D,
 $$
the corresponding standard weight in the unit disc $\D$. The associated weighted Laplace operators  are given by 
$$
L_\alpha= \partial_z \omega_{\p}^{-1} \partial_{\overline{z}} \quad \mbox{ and } L_\alpha^*= 
 \partial_{\overline{z}}   \omega_{\p}^{-1} \partial_z \mbox{ in } \D.
$$

One can check that 

$$L_\p = \p (1-|z|^2)^{-\p-1} \overline{z}  \partial_{\overline{z}} +(1-|z|^2)^{\alpha}\partial_{z}\partial_{\overline{z}}.$$

Let 
$$
T_\alpha:= -\frac{\p^2}{4} \omega_{\p+1}^{-1}+ \frac{1}{2}( L_{\p}+ L^*_{\p}),$$
that is,
$$
T_\p=-\frac{\p^2}{4} (1-|z|^2)^{-\p-1}+\frac{\p}{2}(1-|z|^2)^{-\p-1}\overline{z}\partial_{\overline{z}}+\frac{\p}{2}(1-|z|^2)^{-\p-1}z\partial_z+(1-|z|^2)^{-\p} \partial_{z}\partial_{\overline{z}}.
$$
Thus
$$
T_\p= (1-|z|^2)^{-\p-2} \Delta_{\p/2,\p/2}.
$$

For $\p>-1$, we denote by
\be
M_\p:=M_{\p/2,\p/2}=P_{\p/2,\p/2}[1].\nonumber
\ee
\vskip 1em
Recall that a function $u$ is $T_\p$-harmonic if $T_\p u=0$. Remark that $T_\p$-harmonic functions are exactly 
$(\p/2,\p/2)$-harmonic functions. Thus combining Theorem \ref{2.2} and Theorem \ref{2.3}, we obtain  

 \begin{cor}
Let $\alpha>-1$  and  $ f$ be an  absolutely
continuous such that $\dot{f}\in L^p(\T)$ with $1\leq p\leq \infty$ and  $u=P_{\p/2,\p/2}[f]$. Then 
\begin{enumerate}
    \item[(i)] If $\p>0$, then $\partial_z u$ and $\partial_{\overline{z}} u$ are in $H_G^p(\D)$.
   
    \item[(ii)] If  $\p <0$, then  $\partial_z u$ or $\partial_{\overline{z}}  u$ in $H_G^1(\D)$, if and only if,  $u$ is  zero.
     \item[(iii)] If $\p <0$ and $p<-\frac{1}{\p}$, then $\partial_z u$ and $\partial_{\overline{z}} u$ are in $L^p(\D)$.
     \item[(iv)] If $\p<0$ and $p\geq -\frac{1}{\p}$, then  $\partial_z M_{\p}$   and $\partial_{\overline{z}} M_\p$ $\not\in L^p(\D)$.
\end{enumerate}
\end{cor}
Our results simplify those obtained by Chen et al. in \cite{CCZ}.

\end{document}